\documentclass[11pt]{amsart}
\usepackage{latexsym,amssymb,graphicx,epsfig,amsmath}

\theoremstyle{plain}

\newtheorem{thm}{Theorem}
\newtheorem{prop}[thm]{Proposition}

\newtheorem{corollary}[thm]{Corollary}

\theoremstyle{definition}
\newtheorem{defn}[thm]{Definition}

\begin{document}

\def\Cal#1{{\cal#1}}
\def\<{\langle}\def\>{\rangle}
\def\what{\widehat}\def\wtil{\widetilde}
\def\Z{{\mathbb Z}}\def\N{{\mathbb N}} \def\C{{\mathbb C}}
\def\Q{{\mathbb Q}}\def\R{{\mathbb R}} \def\E{{\mathbb E}}
\def\H{{\mathbb H}}\def\S{{\mathbb S}} 

\def\noproof{\hfill$\square$\break}
\def\noi{\noindent}

\def\Notation{\paragraph{\it Notation.}}
\def\Ackn{\medskip\paragraph{\bf Acknowledgement.}}

\def\Aut{\textsl{Aut}}
\def\Out{\textsl{Out}}
\def\Inn{\textsl{Inn}}
\def\Comm{\textsl{Comm}}
\def\Mod{\textsl{Mod}}
\def\Sym{\textsl{Sym}}
\def\Tv{\textsl{Tv}}

\def\al{\alpha}                 \def\be{\beta}		\def\lam{\lambda}
\def\ga{\gamma}                 \def\Ga{\Gamma}
\def\sig{\sigma}
\def\ep{\epsilon}               \def\varep{\varepsilon}
%

\title{{\bf  Automorphism groups
of some affine and finite type Artin groups}}

\author[R.~Charney]{Ruth Charney}
      \address{Department of Mathematics\\
               Brandeis University\\
               Mail Stop 050\\
               Waltham, MA 02454}
      \email{charney@brandeis.edu}
      \thanks{Charney was partially
      supported by NSF grant DMS-0104026.}

\author[J.~Crisp]{John Crisp}
	\address{I.M.B.(UMR 5584 du CNRS)\\ 
              Universit\'e de Bourgogne\\
		     B.P. 47 870\\
		21078 Dijon, France}
	\email{jcrisp@u-bourgogne.fr}
	
\keywords{}\date{\today}

\date{\emph{August 26, 2004}}

\begin{abstract}
We observe that, for fixed $n\geq 3$, each of the Artin groups of
finite type $A_n$, $B_n=C_n$, and affine type $\wtil A_{n-1}$
and $\wtil C_{n-1}$ is a central extension of a finite index subgroup
of the mapping class group of the $(n+2)$-punctured sphere. (The centre 
is trivial in the affine case and infinite cyclic in the finite type cases).
Using results of Ivanov and Korkmaz on abstract commensurators
of surface mapping class groups we are able to determine the
automorphism groups of each member of these four 
infinite families of Artin groups.
\end{abstract} 

\maketitle

A rank $n$ \emph{Coxeter matrix} is a symmetric $n\times n$ matrix
$M$ with integer entries $m_{ij}\in \N\cup\{\infty\}$ where $m_{ij}\geq 2$
for $i\neq j$, and $m_{ii}=1$ for all $1\leq i\leq n$. 
Given any rank $n$ Coxeter matrix $M$, the \emph{Artin group} of type $M$ is 
defined by the presentation
\[
A(M)\cong \<\ s_1,\dots,s_n\ \mid\ 
\underbrace{s_is_js_i\dots}_{m_{ij}}=\underbrace{s_js_is_j\dots}_{m_{ij}}\ 
\text{ for all } i\neq j, m_{ij}\neq\infty\>\,.
\]
Adding the relations $s_i^2=1$ to this presentation yields a presentation of the
Coxeter group of type $M$ generated by standard reflections $s_i$ and such that
the rotation $s_is_j$ has order $m_{ij}$, for all $1\leq i,j\leq n$. 
A Coxeter matrix $M$ and its Artin group $A(M)$ are said to be of \emph{finite type}
if the associated Coxeter group $W(M)$ is finite, and of 
\emph{affine (or Euclidean) type} if $W(M)$ acts as a proper, cocompact group of
isometries on some Euclidean space with the generators $s_1, \dots , s_n$ acting as 
affine reflections.

The information
contained in the Coxeter matrix $M$ is often displayed in the form of a graph, the 
\emph{Coxeter graph}, whose vertices are numbered $1,..,n$ and which has an edge 
labelled $m_{ij}$ between the vertices $i$ and $j$ whenever $m_{ij}\geq 3$ or $\infty$.
With this particular convention, one usually suppresses the labels which are equal 
to $3$ (but not the corresponding edges!). Note that the absence of an edge between
two vertices indicates that the corresponding generators of $A(M)$ commute.
We say that an Artin group is \emph{irreducible} if its Coxeter graph is connected and
observe that every Artin group is isomorphic to a direct product of irreducible 
Artin groups corresponding to the connected components of its Coxeter graph.

\begin{figure}[ht]
\begin{center}
\includegraphics[width=12cm]{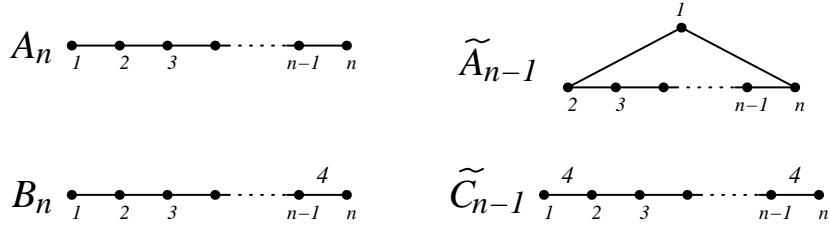}
\end{center}
\caption{Coxeter graphs of finite type $A_n$, $B_n$ and affine type 
$\wtil A_{n-1}$, $\wtil C_{n-1}$, with rank $n\geq 3$. Note that, by convention,
the subscript in each case indicates the geometric dimension of the Artin group.}
\label{Fig1}
\end{figure}

In this paper we concern ourselves with the following four infinite families 
of Artin groups $A=A(M)$ (see Figure \ref{Fig1}):
the finite types $M=A_n$, $B_n$ and affine types $\wtil A_{n-1}$ and $\wtil C_{n-1}$,
of rank $n\geq 3$ in each case. (We refer to \cite{Bou} for the classification
of irreducible finite and affine type Coxeter systems.) For each  Artin group 
$A$ on this list, we determine its automorphism group $\Aut(A)$, its outer automorphism 
group $\Out(A)=\Aut(A)/\Inn(A)$, and the abstract commensurator group $\Comm(A/Z)$
of the group modulo its centre (see Definition \ref{DefnCommG}).

The Artin group $A(A_n)$ is well-known as the braid group on $n+1$ strings, and 
the automorphism groups of the braid groups were determined 
by Dyer and Grossman in \cite{DyGr}.
In this paper we exploit the fact that the Artin groups of 
type $B$, $\wtil A$ and $\wtil C$ also have descriptions as braid 
groups (see \cite{All}) in order to study their automorphism groups. 
Note that the centre of any irreducible Artin group of finite
type is an infinite cyclic group, which we generally denote by $Z$.
(On the other hand, the groups $A(\wtil A_{n-1})$ and $A(\wtil C_{n-1})$ 
have trivial centres.)  As explained in Section \ref{Sect1}, each of the groups 
$A(A_n)/Z$, $A(B_n)/Z$, $A(\wtil A_{n-1})$ and $A(\wtil C_{n-1})$ 
is a finite index subgroup of the mapping class group $\Mod(S_{n+2})$
of the $(n+2)$-punctured 2-sphere $S_{n+2}$. Using this fact, and appealing 
to the work of Ivanov and Korkmaz on the structure of surface mapping 
class groups  \cite{Iv,Kor}, we show:

\begin{thm}\label{main}
For each $n\geq 3$, we have
\begin{itemize}
\item[(i)] $\Comm(G)\cong \Mod(S_{n+2})\ $ for each 
$\ G = A(A_n)/Z,\ A(B_n)/Z,\ A(\wtil A_{n-1}),\ A(\wtil C_{n-1})$.\\
\item[(ii)] 
\begin{itemize}
\item[(a)] \emph{(Dyer, Grossman \cite{DyGr})} 
$\ \Out(A(A_n))\cong\Out(A(A_n)/Z)\cong C_2$
\smallskip
\item[(b)] $\Out(A(B_n))\cong(\Z\rtimes C_2) \times C_2\ $ and
$\ \Out(A(B_n)/Z)\cong C_2\times C_2$
\smallskip
\item[(c)] $\Out(A(\wtil C_{n-1}))\cong \Sym(3)\times C_2$
\smallskip
\item[(d)] $\Out(A(\wtil A_{n-1}))\cong D_{2n}\times C_2\,$.
\smallskip
\end{itemize}
\end{itemize}
\end{thm}

In addition, we show that the exact sequence 
\[ 1\to\Inn(A)\to\Aut(A)\to\Out(A)\to 1\]
 splits
in each of the above cases with the exception of the case $A=A(\wtil C_{n-1})$
with $n\equiv 2 \mod 3$ (see Proposition \ref{nosplit}). In the case 
of $A(\wtil A_{n-1})$, for example, the dihedral factor $D_{2n}$ 
is realised by the group of 
``graph automorphisms'' induced from symmetries of the Coxeter graph, 
which in this case is an $n$-cycle. Every Artin group admits an 
involution $\ep:s_i\mapsto s_i^{-1}$  for all  $i=1,..,n$ which is never inner (because
it reverses the sign of the length function $\ell:A\to\Z$, defined by $s_i\mapsto 1$). 
This accounts for the central $C_2$ factor in each of the outer automorphism groups 
above. 

Note that all of the outer automorphism groups in the theorem are finite with the 
exception of $\Out(A(B_n))$. The infinite cyclic factor in $\Out(A(B_n))$ is generated 
by a ``transvection'' which multiplies each generator $s_i$ by some element of the 
center $Z$. In Section \ref{Sect3}, we study the group of transvections of a finite type 
Artin group. Transvection homomorphisms exist in any Artin group with non-trivial 
center, but in general, they are not automorphisms. In particular, 
transvection automorphisms do not occur when the Artin group has 
abelianisation $\Z$, as in the case of the braid groups. 
However, $A(B_n)$ has abelianisation $\Z\times\Z$ and likewise for Artin groups of 
dihedral type $I_2(m)$ when $m$ is even. (The groups
$\Out(A(I_2(m)))$, for $m\geq 3$, were computed in \cite{GHMR} and are discussed 
in Section \ref{Sect3}). By contrast, transvection homomorphisms appear in the 
abstract commensurator group of \emph{every} finite type Artin group. 
In Proposition \ref{BigComm}, we show that the abstract commensurator 
group of any finite type Artin group contains an infinitely generated free abelian 
group  generated by transvections. 

\Ackn 
The first author would like to thank the Institut de Mathematiques de Bourgogne for 
their hospitality during the development of this paper.

\section{Abstract commensurators and the mapping
class group of a punctured sphere}\label{Sect1}

The following is a special case of Theorem 8.5A of Ivanov's survey article \cite{Iv}.
It is a consequence of ideas laid out in the papers of Ivanov on the mapping class
groups of higher genus orientable surfaces and the work of Korkmaz \cite{Kor}
on the complex
of curves associated to an $m$-punctured sphere. For simplicity we shall write
$S_m$ for the 2-sphere with $m$ points removed. If $S$ is an orientable surface 
(without boundary), we
denote by $\Mod(S)$ the group of all diffeomorphisms of $S$
(not necessarily respecting orientation) modulo diffeotopy.

\begin{thm}[Ivanov, Korkmaz]\label{IvanovKthm}
Let $m\geq 5$. If $\varphi:H\to K$ is an isomorphism between
 finite index subgroups $H,K<\Mod(S_m)$, then $\varphi$ is the restriction to $H$ of 
an inner automorphism of
$\Mod(S_m)$ (conjugation by some element $g\in\Mod(S_m)$).
\end{thm}

\begin{defn}\label{DefnCommG}
Let $G$ be a group. We define the \emph{abstract commensurator group of $G$}
to be 
\[
\Comm(G)=\{ \varphi:H\buildrel\cong\over\rightarrow K\ :\ H,K<G 
\text{ finite index }\}\ /\ \sim
\]
where $\varphi\sim\psi$ if they agree on a finite index subgroup of $G$.
We note that the group structure on $\Comm(G)$ is given by composition of
isomorphisms after appropriate restriction of their domains to finite
index subgroups. 
\end{defn}

\begin{corollary}\label{normalizers}
Let $\Ga=\Mod(S_m)$, $m\geq 5$, let $H<\Ga$ denote a finite index 
subgroup and $N_\Ga(H)$ the normalizer of $H$ in $\Ga$. Then
\begin{itemize}
\item[(i)] $\Comm(H)\cong \Ga$, and
\item[(ii)] $\Aut(H)\cong N_\Ga(H)$.
\end{itemize}
\end{corollary}

\begin{proof}
(i) There is a natural homomorphism $\Ga \to \Comm(H)$ which takes $\gamma \in \Ga$ to 
conjugation by $\gamma$ restricted to the finite index subgroup 
$H \cap \gamma^{-1}H\gamma$. It follows from Theorem \ref{IvanovKthm} that this map is 
surjective. To show that it is injective, we observe that $\Ga$ has trivial 
``virtual centre'', that is, any element  which centralizes a finite index subgroup
is necessarily trivial. This is because any mapping class $g\in\Mod(S)$
is determined by its action on a finite number of isotopy classes of
simple closed curves on the surface $S$, and therefore by its action
on the $n$-th powers $T_i^n$ of the corresponding Dehn twists for
any choice of $n$.
Given a finite index subgroup $K<\Ga$ centralized by $g$ 
we can always find a sufficiently large $n$ such that $K$ contains all
$T_i^n$. Then $g$ is determined by its action on $K$ and so must be the
identity in $\Ga$.

(ii) There is now an obvious homomorphism $\Aut(H)\to\Ga$ which factors
through the inverse of the isomorphism just described, and realizes each 
automorphism of $H$ as conjugation by an element of $\Ga$. 
It follows easily that this map is injective with image precisely the 
normalizer $N_\Ga(H)$ of $H$ in $\Ga$.
\end{proof}

\section{Some interesting finite index subgroups of $\Mod(S_m)$
and their automorphisms}\label{Sect2}

Let $\Sym(m)$ denote the group of permutations of the set $\{1,2,..,m\}$, for $m\geq 1$.
If $k<m$ then we consider $\Sym(k)$ as a subgroup of $\Sym(m)$ consisting of
those permutations which fix the subset $\{k+1,..,m\}$. 
Notice that each mapping class of the $m$-punctured sphere $S_m$
induces a permutation of the punctures. Moreover, the mapping class group acts
on the orientation class of the surface. Thus we have a surjective 
group homomorphism $\pi:\Mod(S_m)\to\Sym(m)\times C_2$ with kernel the group of
pure orientation preserving mapping classes. Note that, strictly speaking, the 
definition of $\pi$ depends on a labelling of the punctures by the numbers $1$ 
through $m$.  

\medskip

\noindent{\bf Notation.}
Fix $n\geq 3$. We shall write $\Ga=\Mod(S_{n+2})$ and $\Ga_A$, $\Ga_B$, 
$\Ga_{\wtil C}$ respectively, for the finite index subgroup of $\Ga$ generated
by those orientation preserving diffeomorphisms which fix the last
1, 2, 3 punctures of $S_{n+2}$ respectively. These three subgroups are, in other words,
the preimages under $\pi$ of the subgroups $\Sym(k)<\Sym(n+2)\times C_2$,
with $k=n+1,n,n-1$ respectively. 
In particular we have $[\Ga :\Ga_A]=2(n+2)$, $[\Ga :\Ga_B]=2(n+1)(n+2)$
and $[\Ga :\Ga_{\wtil C}]=2n(n+1)(n+2)$.

\bigskip
  
It is well known that the $(n+1)$-string braid group, $A(A_n)$, is isomorphic to 
the mapping class group of the $(n+1)$-punctured disk  \emph{relative} to the boundary 
of the disk. (That is, diffeotopies are required in this case to fix the boundary 
of the disk \emph{pointwise}). The centre $Z$ of this mapping class group is 
generated by the Dehn twist about a curve parallel to the boundary, and is precisely 
the kernel of the natural homomorphism to the mapping class group $\Ga=\Mod(S_{n+2})$ 
(induced by inclusion of the punctured disk in $S_{n+2}$). Thus we may realise the 
group $A(A_n)/Z$ as the finite index subgroup $\Ga_A$ of $\Ga$.

It is less well-known (see \cite{All}) that the Artin groups $A(B_n)$ and
$A(\wtil C_{n-1})$, respectively, are isomorphic  to the subgroups of the braid
group $A(A_n)$ leaving fixed one, respectively two, of the punctures in the disk.
We note that the latter of these two groups intersects the centre of the
braid group trivially, from which we deduce that
\begin{itemize}
\item [$\bullet$] $\Ga_A \cong A(A_n)/Z$,
\item [$\bullet$] $\Ga_B \cong A(B_n)/Z$, and
\item [$\bullet$] $\Ga_{\wtil C} \cong A(\wtil C_{n-1})$.
\end{itemize}

\medskip

\begin{figure}[ht]
\begin{center}
\includegraphics[width=8cm]{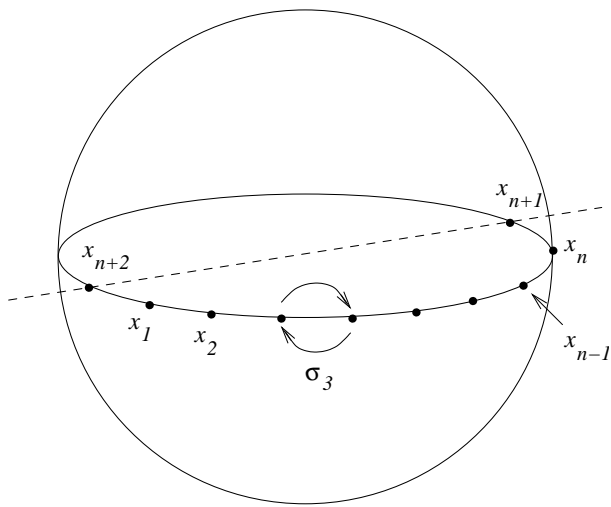}
\end{center}
\caption{Finite index subgroups of $\Ga=\Mod(S_{n+2})$ isomorphic to $A(A_n)/Z$, 
$A(B_n)/Z$ and $A(\wtil C_{n-1})$, $n\geq 3$, fix the subsets $\{x_{n+2}\}, \{x_{n+1},x_{n+2}\}$ and $\{x_1,x_{n+1},x_{n+2}\}$ respectively.}\label{Fig2}
\end{figure}

These isomorphisms may be described explicitly as follows. Distribute the punctures 
$x_i$ ($i=1,..,n+2$) evenly along half the equator in the order 
$x_{n+2},x_1,x_2,..,x_n,x_{n+1}$ so that $x_{n+1}$ and $x_{n+2}$ are antipodal 
(see Figure \ref{Fig2}). Then the $i$th standard generator of $A(A_n)$ maps to the
positive braid twist $\sigma_i$ exchanging $x_i$ with $x_{i+1}$. 
The $i$th standard generator of $A(B_n)$ maps to $\sigma_n^2$ if $i=n$
and $\sigma_i$ otherwise. The $i$th standard generator of $A(\wtil C_{n-1})$ maps 
to $\sigma_i^2$ if $i\in\{ 1,n\}$  and $\sigma_i$ if $i\in\{ 2,..,n-1\}$. (This maps  
$A(\wtil C_{n-1})$ into the subgroup of $\Ga$ fixing $x_{n+2}, x_{n+1}$, and $x_1$. 
Renumbering the punctures gives  $\Ga_{\wtil C}$.)

\begin{prop}\label{ABC}
\begin{description}
\item[(i)] $\Comm(\Ga_A)\cong\Comm(\Ga_B)\cong\Comm(\Ga_{\wtil C})\cong \Ga\,,$
\item[(ii)] there exist short exact sequences 
\[
\begin{aligned}
1\to \Ga_A \to &\Aut(\Ga_A)\to C_2\to 1\\
1\to\Ga_B \to &\Aut(\Ga_B)\to C_2\times C_2\to 1\\
1\to\Ga_{\wtil C} \to &\Aut(\Ga_{\wtil C})\to \Sym(3)\times C_2\to 1
\end{aligned}
\]
\end{description}
\end{prop}
 
\begin{proof}
Let $H$ denote one of $\Ga_A$, $\Ga_B$, or $\Ga_{\wtil C}$, and $P$ the set of 
1,2 or 3 punctures fixed by all elements of $H$.
It is easily checked that $N_\Ga(H)$ is exactly the group of all mapping
classes (orientation preserving or not) which leave $P$ setwise invariant.
It follows that $N_\Ga(H)/H$ is a direct product of $C_2=\Ga/\Ga^+$ with the
group $\Sym(P)$ of permutations of the set $P$. Statements (i) and (ii)
now follow immediately from the corresponding statements in
Corollary \ref{normalizers}.   
\end{proof}

It is easy to see, with the aid of Figure \ref{Fig2},
that the short exact sequences in the proposition above are split 
in the cases $\Ga_A$ and $\Ga_B$. The group $C_2\times C_2$ may be generated
by a reflection in the equatorial plane and a rotation about an axis in the 
equatorial plane which exchanges the points $x_{n+1}$ and $x_{n+2}$. Only
the reflection leaves $\Ga_A$ invariant. 
 
For $\Ga_{\wtil C}$, the exact sequence splits 
providing $n+2 \equiv 0,2$ mod $3$. To see this, think of $S_{n+2}$ as a sphere 
with $3k$ punctures arranged symmetrically along the equator and, if $n+2 \equiv 2(3)$, 
the remaining two punctures at the north and south poles. Then there is an 
orientation-preserving action of $\Sym(3)$ on $S_{n+2}$ generated by an order 
three rotation $\rho$ about the vertical axis through the poles and an order two 
rotation $\eta$ about a horizontal axis through one ($k$ odd) or two ($k$ even) 
punctures at the equator.  If $x$ is a puncture fixed by $\eta$, then the orbit 
of $x$ under $\rho$ is preserved by $\Sym(3)$, hence for an appropriate numbering
of the punctures (namely, so that $P=\{ x,\rho(x),\rho^2(x)\}$)
the action of $\Sym(3)$ lies in $\Aut(\Ga_{\wtil C})$ and is faithful
on $P$. This splits the $\Sym(3)$ factor. The $C_2$ factor is realized by reflection 
through the equatorial plane.

\begin{prop}\label{nosplit}
If $\Ga_{\wtil C}=A(\wtil C_{n-1})$ with $n+2\equiv 1$ mod $3$, 
then $\Aut(\Ga_{\wtil C})$ does not contain a subgroup isomorphic to $\Sym(3)$. 
In this case, the exact sequence for $\Ga_{\wtil C}$ in the previous proposition
does not split. 
\end{prop}

\begin{proof}
Suppose  $\Aut(\Ga_{\wtil C}) \subset \Ga$ contains a subgroup $K$ 
isomorphic to $\Sym(3)$ and let $H$ be the subgroup of $K$ generated by a 3-cycle. 
Then $H$ is a subgroup of $\Mod(S_{n+2})$, hence $H$ acts by permutations on the set 
of punctures $\{x_1, \dots x_{n+2}\}$. Every $H$ orbit consists of either $1$ or $3$ 
punctures. Since $n+2 \equiv 1 (3)$, the number of fixed points of $H$ must also be 
congruent to $1$ mod $3$.  If $H$ fixes $4$ or more punctures, then it lies in a 
subgroup of $\Ga$ isomorphic to $\Ga_{\wtil C}$ which is impossible since  
$\Ga_{\wtil C}$ is torsion-free (since, by \cite{All}, it is a subgroup of 
the braid group on $n+1$ strings which is known to be torsion free).

Thus $H$ must have a unique fixed puncture. Since $H$ is normal in $K$, all of $K$ 
must fix this puncture and hence $K$ lies in a subgroup of $\Ga$ isomorphic to $\Ga_A$. 
By a theorem of Bestvina, \cite{Best} Thm 4.5, every finite subgroup of $\Ga_A$ 
is cyclic, so we arrive at a contradiction.
\end{proof}

\begin{figure}[ht]
\begin{center}
\includegraphics[width=8cm]{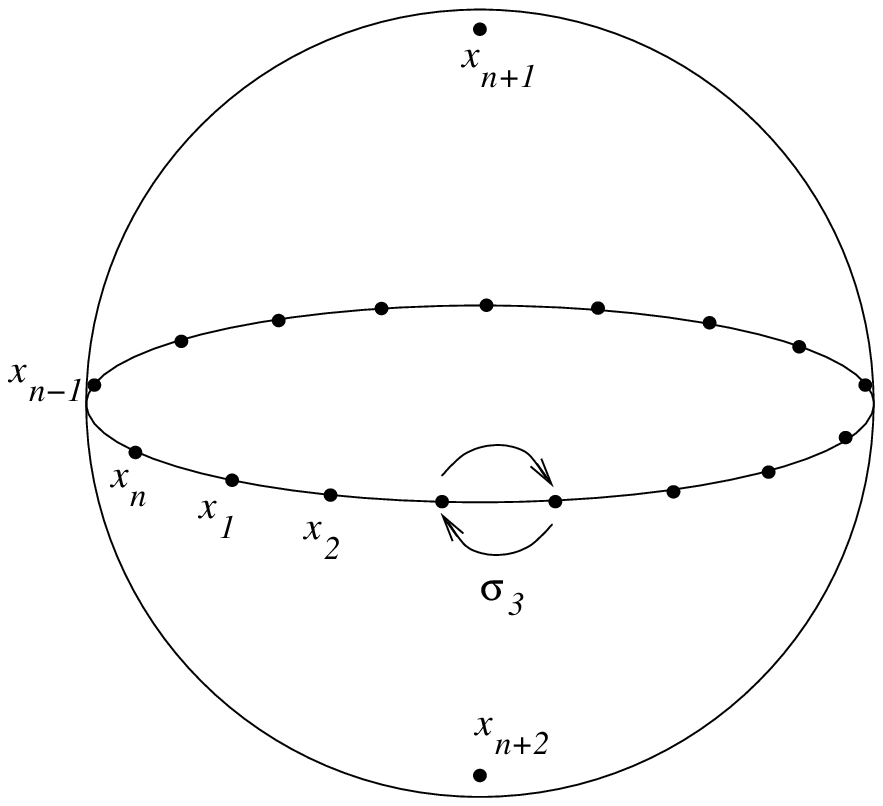}
\end{center}
\caption{The finite index subgroup $\Ga_{\wtil A}$ of $\Ga=\Mod(S_{n+2})$ 
isomorphic to $A(\wtil A_{n-1})$.}\label{Fig3}
\end{figure}

It has been observed by several authors \cite{All}\cite{tD}\cite{KP} that there exists a semidirect product decomposition 
\[
A(B_n)\cong A(\wtil A_{n-1})\rtimes\Z\,.
\]
(See \cite{ChPe} for further discussion of this decomposition.)
The Coxeter graph of $A(\wtil A_{n-1})$ is an $n$-cycle and the generator of the 
cyclic factor in the semi-direct product acts on $A(\wtil A_{n-1})$ via 
an order $n$ rotation of this graph. The centre of $A(B_n)$ is the subgroup
$n\Z$ of the cyclic factor. Thus we also have
\[
\Ga_B\cong A(\wtil A_{n-1})\rtimes\Z/n\Z\,.
\]
Here we interpret $A(\wtil A_{n-1})$ as the subgroup of $\Ga_B$ consisting of
mapping classes of ``zero angular momentum'' about the axis through the two fixed
punctures.  More precisely, suppose 
that the punctures of $S_{n+2}$ are arranged so that $x_{n+1}$ and $x_{n+2}$ are
placed at the north and south pole respectively, and the remaining points $x_1,..,x_n$
are equally spaced (in that cyclic order) around the equator, see figure \ref{Fig3}. 
Define $\Ga_{\wtil A}$
to be the subgroup of $\Ga$ generated by the braid twists $\sigma_i$ which exchange
the points $x_i$ and $x_{i+1}$, for $i=1,..,n$, indices taken mod $n$. 
The isomorphism $A(\wtil A_{n-1})\cong\Ga_{\wtil A}$ is given by sending the $i$th
standard generator of $A(\wtil A_{n-1})$ to $\sigma_i$. 
 Note that $\Ga_{\wtil A}$ is a finite index subgroup of $\Ga$ of 
index $[\Ga:\Ga_{\wtil A}]=n[\Ga:\Ga_B]=2n(n+1)(n+2)$.
In particular, $\Comm(\Ga_{\wtil A})\cong \Ga$.

\begin{prop}\label{Atilde}
Let $\Ga_{\wtil A}\cong A(\wtil A_{n-1})$ be as above. 
Then $N_\Ga(\Ga_{\wtil A})=N_\Ga(\Ga_B)$. Consequently
$\Aut(\Ga_{\wtil A})\cong\Aut(\Ga_B)$ and there is an
exact sequence
\[
1\to\Ga_{\wtil A} \to \Aut(\Ga_{\wtil A})\to D_{2n}\times C_2\to 1\,,
\]
where $D_{2n}$ denotes the dihedral group of order $2n$.
Moreover, a splitting of this sequence is given by realizing the group of
outer automorphisms as the direct product of the group of graph automorphisms
of $A(\wtil A_{n-1})$ with the group of order $2$ generated by the inversion
$\ep$ which inverts each standard generator.
\end{prop}

\begin{proof}
We suppose that the punctures of $S_{n+2}$ are arranged as shown
in Figure \ref{Fig3}. The outer automorphism group of $\Ga_B$ may 
then be generated by the reflection in the equatorial plane and a rotation
exchanging the two poles (this gives a  splitting of the exact
sequence of Proposition \ref{ABC} slightly different to the one previously mentioned).
It follows easily that $\Ga_{\wtil A}$ is invariant by these mapping classes, 
and hence a characteristic
subgroup of $\Ga_B$. In particular, we have $N_\Ga(\Ga_B)< N_\Ga(\Ga_{\wtil A})$.

Let $P=\{x_{n+1},x_{n+2}\}$. To prove that $N_\Ga(\Ga_{\wtil A})=N_\Ga(\Ga_B)$ 
it suffices now to check that every mapping class which normalizes $\Ga_{\wtil A}$
leaves $P$ invariant. But if $g\in N_\Ga(\Ga_{\wtil A})$ then 
$\pi(g)$ normalizes $\pi(\Ga_{\wtil A})=\Sym(n)\times C_2$, and
this is just to say
that $g$ induces a symmetry of the puncture set which leaves $P$
invariant, as required.

The remaining statements are now easily deduced from Corollary
\ref{normalizers} and Proposition \ref{ABC}.
\end{proof}

In summary, Propositions \ref{ABC} and \ref{Atilde} give us explicitly the
automorphism groups of $A(A_n)/Z$, $A(B_n)/Z$, $A(\wtil A_{n-1})$, and
$A(\wtil C_{n-1})$, for all $n\geq 3$, as well as the fact that the
abstract commensurator group is in each case isomorphic to $\Mod(S_{n+2})$.
This proves Theorem \ref{main} with the exception of the
first isomorphism in each of part (ii)(a) and (b). 
In the next section we shall recover the group $\Aut(A(B_n))$ from the
above result on $\Aut(A(B_n)/Z)$. In a similar fashion, the result of \cite{DyGr}
on automorphisms of the braid group is easily recovered from the computation of
$\Aut(A(A_n)/Z)$. These steps will complete the proof of Theorem \ref{main}. 

\section{Transvections and automorphisms of $A(B_n)$}\label{Sect3}

Let $H$ be a group with nontrivial centre $Z$. By a \emph{transvection} of $H$
we mean a homomorphism $T_\lam:H\to H$ given by $T_\lam(x)=x\lam(x)$,
where $\lam$ denotes a function $H\to Z$. 
The fact that $T_\lam$ is a homomorphism requires that $\lam$ is 
also a homomorphism, and it is easily checked that the composite of two 
transvections is again a transvection.  We note, however, that transvections are
 typically not automorphisms. It is easily
shown that a transvection is an automorphism (resp. injective,
resp. surjective) if and only if its restriction to the centre $Z$ is 
an automorphism (resp. injective, resp. surjective). 

Let $\Tv(H)$ denote the group of automorphisms 
of $H$ which are transvections. Then, almost by definition, we have
an exact sequence
\[
1\to \Tv(H)\to \Aut(H)\to \Aut(H/Z)\,.
\]
Note that $\Inn(H)\cong H/Z$ always maps isomorphically onto $\Inn(H/Z)$.
In general there is an issue, however, as to whether \emph{every} 
outer automorphism of $H/Z$ lifts to an automorphism of $H$, 
and as to whether the extension splits.

In the case of an irreducible finite type Artin group $A$, the extension 
splits at least over a finite index subgroup $\Aut^*(A/Z)$ of $\Aut(A/Z)$ 
defined as follows. Let $\ell : A \to \Z $ be the length homomorphism which 
maps each generator of $A$ to $1$. Recall that the  centre $Z$ of $A$ is 
infinite cyclic. Let $\zeta$ be a generator of $Z$ and let $d=\ell(\zeta)$. 
The length homomorphism on $A$ descends to a ``length homomorphism''
$\bar{\ell} : A/Z \to \Z/d\Z $. Define $\Aut^*(A)$, respectively $\Aut^*(A/Z)$, 
to be the group of length preserving or reversing automorphisms, that is, 
automorphisms $\alpha$ such that $\ell\circ\alpha = \pm\ell$ 
(respectively $\bar\ell\circ\alpha = \pm\bar\ell$).

\begin{prop}\label{autstar}
Suppose $A$ is an irreducible Artin group of finite type. Then the natural map 
$\pi : \Aut^*(A) \to \Aut^*(A/Z)$ is an isomorphism.  If, in addition, the 
abelianization of $A$ is infinite cyclic (i.e., $A$ is not of type 
$I_2(2m), B_n$, or $F_4$), then $\Tv(A)$ is trivial and $\Aut(A) \cong \Aut^*(A/Z)$.
\end{prop}

\begin{proof}
Note that a non-trivial transvection $T_\lambda$ of $A$ can never be length preserving, 
and can only be length reversing if $d=2$ (and $\lambda(\cdot)=\zeta^{-\ell(\cdot)}$). 
Since $d>2$ for all
irreducible finite type Artin groups we conclude that $\ker(\pi)=\Tv(A)\cap\Aut^*(A)$
is trivial. 
 
For any length preserving $\alpha \in \Aut^*(A/Z)$, we define a lift 
$\hat\alpha \in \Aut^*(A)$ as follows.  (The lift for length reversing 
automorphisms is defined analogously). By definition,
for any $a \in A$, $\alpha (aZ)=bZ$ for some $b \in A$ with 
$\ell(b) \equiv \ell(a)$ mod $d$. Since $Z$ is generated by an element of 
length $d$, the coset $bZ$ contains a unique representative $b'$ satisfying 
$\ell(b') = \ell(a)$. Define $\hat\alpha (a)=b'$. To verify that $\hat\alpha$ 
is a homomorphism, note that if $\hat\alpha(a_1)=b_1$ and $\hat\alpha(a_2)=b_2$, then
\[
\begin{aligned}
\alpha(a_1a_2Z)&=\alpha(a_1Z)\cdot \alpha(a_2Z)=b_1Z\cdot b_2Z=b_1b_2Z, \\
\ell(a_1a_2)&=\ell(a_1)+\ell(a_2)=\ell(b_1)+\ell(b_2)=\ell(b_1b_2),
\end{aligned}
\]
so it follows that $\hat\alpha(a_1a_2)=b_1b_2$. 
Moreover, it is straightforward to verify that if $\alpha_1$ and $\alpha_2$ 
are two elements of $\Aut^*(A/Z)$, then $\alpha_1\circ\alpha_2$ lifts to 
$\hat\alpha_1\circ\hat\alpha_2$. In particular, the lift of an automorphism 
is also an automorphism. Thus, $\alpha \mapsto \hat\alpha$ defines a section 
for the projection $\pi$, and $\pi$ is an isomorphism.

If the abelianization of $A$ is infinite cyclic, then the length homomorphism 
$\ell$ can be identified with the abelianization map $A \to A^{\text{ab}}\cong \Z$. 
Hence every automorphism of $A$ is either length preserving or length reversing. 
The second statement of the proposition follows.
\end{proof}

In the case of $A=A(A_n)$, it follows from Proposition \ref{ABC} that 
$\Aut^*(A/Z)=\Aut(A/Z) \cong \Ga_A \rtimes C_2$ which, combined with 
Proposition \ref{autstar}, yields the main result of  \cite{DyGr}, 
that $\Out(A(A_n))=C_2$. 

 We remark, however, that while $\Aut(A)=\Aut^*(A)\cong \Aut^*(A/Z))$ 
 when $A$ has infinte cyclic abelianization, it is not always the case that
$\Aut(A)\cong \Aut(A/Z)$ (or equivalently that $\Out(A) \cong \Out (A/Z)$). 
 In particular, for $A$ of type $I_2(m)$ with $m$ odd,
 $A/Z\cong C_2\star C_m$, and $\Out(A/Z)$ is isomorphic
to the group of units of $\Z/m\Z$, while $\Out(A)=C_2$.

In the case of an Artin group whose abelianization is not infinite cyclic, 
the situation is more complicated.

\begin{prop}\label{trans}
If $A$ is an irreducible Artin group of type $I_2(2n), B_n$, $n\geq 3$, or $F_4$, 
then $\Tv(A) \cong \Z$. If $A$ is of type $I_2(4)\equiv B_2$ 
then $\Tv(A)\cong D_\infty$. Moreover, in all these cases $Aut(A)$ 
contains $\Tv(A) \rtimes Aut^*(A)$ as a finite index subgroup.
\end{prop}

\begin{proof}
We refer the reader to \cite{Bou} and \cite{BS} for descriptions of the 
Coxeter graphs, and central elements respectively, for each of the irreducible
finite type Artin groups. In each case listed above, the abelianization of $A$ 
is $\Z \times \Z$, and the generator $\zeta$ of $Z$ maps to an element
$(r,s)\in \Z\times\Z$, where $(r,s)=(n,n)$ for type $I_2(2n)$, $(12,12)$ for type
$F_4$, and $(n(n-1),n)$ for type $B_n$. Any homomorphism $\lam:A\to Z$ 
can be obtained by composing the abelianisation homomorphism with 
a map $\Z\times\Z\to Z$ given by $(r,s)\mapsto \zeta^{pr+qs}$ for a pair 
of integers $p,q\in\Z$. The associated transvection then satisfies 
$T_\lam(\zeta)=\zeta^{k}$ where $k=1+(p+q)n$, or $1+(p(n-1)+q)n$,
depending on the case ($n=12$ for type $F_4$).
On the other hand $T_\lam$ is an automorphism if and only if $k=\pm 1$. 
When $n\geq 3$ this is only possible if $k=1$ and $p+q=0$, or 
$(p(n-1)+q)=0$, respectively.  
This allows exactly an infinite cyclic group of transvection 
automorphisms (all of which act by the identity on $Z$). 

In the remaining case, type $I_2(4)\equiv B_2$, we have the constraint
$k=1+2(p+q)=\pm 1$, and so $p+q=0$ or $-1$. We note that $A=A(I_2(4))$ is given 
by the presentation $\<a,b\mid abab=baba\>$, where $\zeta=abab$. 
A transvection of the second type ($p=0,q=-1$) is realized by the 
automorphism $T_0:a\mapsto a$ and $b\mapsto (aba)^{-1}$.
Since in general we have 
\[
T_\lam\circ T_\mu=\begin{cases}
T_{\lam+\mu} \text{ if } T_\lam(\zeta)=\zeta\\
T_{\lam-\mu} \text{ if } T_\lam(\zeta)=\zeta^{-1}
\end{cases}
\]
it is easily checked that $\Tv(A)$ is an infinite dihedral group with $T_0$ acting
as a direction reversing involution.

The subgroup $\Aut^*(A/Z)$ is finite index in $\Aut(A/Z)$ (since $\bar\ell$
takes values in the finite group $\Z/d\Z$ and $A/Z$, being finitely generated, 
admits only a finite number of homomorphisms to $\Z/d\Z$). Hence the
inverse image $\pi^{-1}(\Aut^*(A/Z)) < \Aut(A)$ is a finite index subgroup 
of $\Aut(A)$. By the previous proposition, this subgroup splits as a semi-direct 
product $Tv(A) \rtimes Aut^*(A)$.
\end{proof}

For $A$ of type $I_2(m)$, we have the following presentation:
\[
A(I_2(m))=\<\ a,b\ \mid\ \underbrace{aba\dots}_m=\underbrace{bab\dots}_m\ \>\,.
\]
In this case the outer automorphism groups were computed in  \cite{GHMR}:
\[
\begin{aligned}
Out(A(I_2(m)) &\cong C_2, &\text{$m$ odd}\\
Out(A(I_2(m)) &\cong (\Z \rtimes \<\gamma\>)\times C_2 \cong D_{\infty}\times C_2\,,
 &\text{$m$ even}
\end{aligned}
\]
where $\gamma$ denotes the graph involution $\gamma: a\leftrightarrow b$, 
and the $C_2$ factor in each case is generated by the inversion automorphism 
$\epsilon:a\mapsto a^{-1},b\mapsto b^{-1}$. In the case $m$ even,
the infinite cyclic factor $\Z$ is generated by the automorphism 
$\eta:a\mapsto aba$ and $b\mapsto a^{-1}=b(ab)^{-1}$, 
and contains the group of transvections as a subgroup of index $m/2$, when
$m\geq 6$. (Note that $\eta$ fixes both elements $ab$ and $ba$, and that the 
centre is generated by $(ab)^{m/2}=(ba)^{m/2}$).
It follows that, in the case $A=A(I_2(m))$ with $m\geq 6$ even, the subgroup
$Tv(A) \rtimes Aut^*(A)$ has index $m/2$ in the entire automorphism group.
This behaviour seems to be repeated in the case $A=A(B_n)$, for $n\geq 3$, 
where, as we will see below, $\pi^{-1}(\Aut^*(A/Z))$  is a subgroup of 
index $2$ in $\Aut(A)$.  The type $I_2(4)\equiv B_2$ Artin group
is however exceptional in this regard. In this case the group of 
transvections is generated by $\eta^2$ and the element $T_0=\ep\circ\eta\circ\ga$,
from which it follows that $Tv(A) \rtimes Aut^*(A)$ is the whole automorphism group.
We do not know what the index is (or whether the subgroup 
is proper) in the case of $A(F_4)$.

\begin{prop}\label{autosB}
Let $n\geq 3$ and $A=A(B_n)$. Then
\[
\begin{aligned}
 Aut(A) &\cong (\Ga_B \times \Tv(A)) \rtimes (C_2 \times C_2)\\
Out(A) &\cong (\Z \rtimes C_2) \times C_2
\end{aligned}
\]
\end{prop} 

\begin{proof}
We first show that the map $\phi:\Aut(A)\to \Aut(A/Z)$ is surjective and splits.
Recall from Section \ref{Sect2} that $A/Z\cong \Ga_B\cong \Ga_{\wtil A}\rtimes\Z/n\Z$. 
It is clear from the proof of Proposition \ref{Atilde} that elements of $\Aut(\Ga_B)$ 
act on this semidirect product by leaving the characteristic subgroup $\Ga_{\wtil A}$
invariant and by mapping the generator of the cyclic factor $\Z/n\Z$ either to itself
or its inverse. Any such automorphism lifts uniquely to an automorphism of 
$A(B_n)\cong \Ga_{\wtil A}\rtimes\Z$. Thus we have a section to the map $\phi$.
 
By Proposition \ref{ABC}, this proves that    
\[
\Aut(A)\cong\Tv(A)\rtimes\Aut(\Ga_B) \cong \Tv(A)\rtimes(\Ga_B\rtimes(C_2 \times C_2)).
\]
A direct check shows that inner automorphisms, $\Ga_B$, commute with transvections, 
and the first statement of the proposition follows. 

We claim that the $C_2$ factor generated by the inversion $\ep$ also commutes with 
transvections. This follows from the fact that $\ep$ restricts to a length reversing
automorphism of the center $Z$ and hence takes the generator $\zeta$ of $Z$ to 
$\zeta^{-1}$. For any transvection $T$ and any generator $s_i$,  $T(s_i)=s_i\zeta^k$ 
for some $k$, so
\[
\ep\circ T \circ \ep (s_i)  =\ep\circ T(s_i^{-1})
=\ep(s_i^{-1}\zeta^{-k})=s_i\zeta^k=T(s_i).
\]
By Proposition \ref{trans}, $\Tv(A) \cong \Z$  and the second claim of the proposition 
now follows.
\end{proof}

\noindent {\bf Remark.}  The $\Z$ factor in the decomposition 
$A(B_n)=\Ga_{\wtil A}\rtimes \Z$ is generated by the element $\delta=s_1s_2\dots s_n$ 
of $A(B_n)$ and it acts on $\Ga_{\wtil A}$ via rotation of the $\wtil A$ graph. These 
rotations thus become inner automorphisms in $Aut(A(B_n))$.  The remaining graph 
automorphisms (the reflections) give rise to one of the $C_2$ factors in $Aut(A(B_n))$.
(The other $C_2$ factor is generated by the inversion $\ep$.) A generator of this $C_2$ 
is represented by the automorphism $\tau$ defined by
\[
\begin{aligned}
\tau(s_i)&=s_{n-i},\,\, i=1,\dots , n-1\\
\tau(\delta)&= \delta^{-1}
\end{aligned}
\]
and hence 
$\tau(s_n)=\tau((s_1\dots s_{n-1})^{-1}\delta)=(\delta s_{n-1}\dots s_1)^{-1}$. 
In particular,  $\tau$ does not lie in $Tv(A)\rtimes \Aut^*(A)$ (as long as $n\geq 3$),
and this latter group is index 2 in the full automorphism group. 
(By contrast, when $A$ is of type $B_2\equiv I_2(4)$, the 
involution $\tau$ just described is exactly the exceptional transvection
$T_0=\ep\circ\eta\circ\gamma$ already mentioned in the preceding discussion).
   
\bigskip

If $H$ is a group whose centre $Z$ is a finite rank free 
abelian group, the injective transvections
which are not automorphisms of $H$ have as image a proper finite
index subgroup of $H$ and therefore represent infinite order 
elements of $\Comm(H)$ (\emph{not} equivalent to automorphisms).
In the case of a finite type Artin group, this implies that the abstract 
commensurator group must be large.

\begin{prop} \label{BigComm}
Let $A$ be any finite type Artin group. Then $\Comm(A)$ contains an infinitely 
generated free abelian subgroup generated by transvections.
\end{prop}

\begin{proof}
It suffices to prove the proposition for irreducible $A$.
Let $\ell : A \to \Z $ be the length homomorphism. Let $\zeta$ be the generator 
of $Z(A)$ of positive length $d=\ell(\zeta)$.  Then for every $m \in \Z$, there 
is a homomorphism $\lam_m: A \to Z(A)$ defined by $\lam_m(a)=\zeta^{m\ell(a)}$.  
The associated transvection takes $a$ to $a\zeta^{m\ell(a)}$. In particular, it 
takes $\zeta$ to $\zeta^{md+1}$. We denote this transvection by $T_{md+1}$.  

It is straightforward to verify that $T_{md+1} \circ T_{kd+1} = T_{(md+1)(kd+1)}$. 
Thus, $T_{md+1}$ and $T_{kd+1}$ are commuting elements of $\Comm(A)$ and to find an 
infinitely generated free abelian subgroup of $\Comm(A)$, it suffices to find an 
infinite sequence of integers $\{n_i\}$ which are relatively prime and satisfy 
$n_i \equiv 1$ mod $d$.  This can be done inductively by setting $n_1=1+d$ and 
$n_{i+1}=1+d(n_1n_2\dots n_i)$. The transvections $T_{n_i}$ are then linearly 
independent.
\end{proof}


\begin{thebibliography}{99}

\bibitem{All}
\textsc{D. Allcock}. Braid pictures for Artin groups. 
\emph{Trans. Amer. Math. Soc.} {\bf 354} (2002),
3455--3474.

\bibitem{Best}
\textsc{M. Bestvina}. Non-positively curved aspects of Artin groups of finite type.
\emph{Geometry \& Topology} {\bf 3} (1999), 269--302.

\bibitem{Bou}
\textsc{N. Bourbaki}. Groupes et alg\`ebres de Lie, Chaps 4--6, Hermann, Paris, 1968.

\bibitem{BS}
\textsc{E. Brieskorn} and \textsc{K. Saito}, Artin-Gruppen und Coxeter-Gruppen.
\emph{Invent. Math.} {\bf 17} (1972),245\hbox{--}271.

\bibitem{ChPe}
\textsc{R. Charney} and \textsc{D. Peifer}. The $K(\pi,1)$ 
conjecture for the affine braid groups.
\emph{Comm. Math. Helv.} {\bf 78} (2003), 584-600.

\bibitem{tD}
\textsc{T. tom Dieck}. Categories of rooted cylinder ribbons and their representations.
\emph{J. Reine Angew. Math.}{\bf 494}(1998), 35--63.

\bibitem{DyGr}
\textsc{J.L. Dyer} and \textsc{E.K. Grossman}.
The automorphism groups of the braid groups.
\emph{Amer. J. Math.} {\bf 103} (1981), 1151--1169.

\bibitem{GHMR}
\textsc{N.D. Gilbert}, \textsc{J. Howie}, \textsc{V. Metaftsis} and \textsc{E. Raptis}.
Tree actions of automorphism groups,
\emph{J. Group Theory} {\bf 3} (2000), no.2, 213--223.

\bibitem{Iv}
\textsc{N.V. Ivanov}. 
Mapping Class Groups, \emph{in} ``Handbook of Geometric Topology" 
\emph{ed.} R.J.~Daverman and R.B.~Sher, Elsevier, 2002.

\bibitem{KP}
\textsc{R. Kent} and \textsc{D. Peifer}. 
A geometric and algebraic description of the annular braid groups.
\emph{Int. Jour. Alg. Comp.}{\bf 12}(2002), 85--97.

\bibitem{Kor}
\textsc{M. Korkmaz}. 
Complexes of curves and mapping class groups.
Ph.D. Thesis, Michigan State University, 1996.
 
\end{thebibliography}
\end{document}